\documentclass[11pt]{article}
\usepackage{mathrsfs}

\newcounter{example}

\usepackage{latexsym}
\usepackage{amsmath}
\usepackage{amssymb}
\usepackage{graphicx}
\usepackage{textcomp}
\usepackage{multirow}
\usepackage{enumerate}
\usepackage{cancel}
\usepackage{hyperref}
\usepackage{caption}  
\usepackage[usenames,dvipsnames]{color}
\usepackage[left, pagewise]{lineno}
\usepackage[ruled,vlined]{algorithm2e}
\usepackage{epstopdf}
\graphicspath{{./Figures/}{./Figures2/}}
\newcommand{\mbr}{\mathbb{R}}
\newcommand{\bP}{{\bf P}}

\newcommand{\bJ}{{\bf J}}

\title{Structure-Preserving Galerkin POD-DEIM Reduced-Order Modeling of Hamiltonian Systems}
\author{Zhu Wang \thanks{ Corresponding author.
           Department of Mathematics, University of South Carolina, Columbia, SC 29208.
           Email: {wangzhu@math.sc.edu}. }
}
\date{}
\begin{document}
\maketitle
\begin{abstract}
A structure preserving proper orthogonal decomposition reduce-order modeling approach has been developed in \cite{gong2017structure} for the Hamiltonian system, which uses the traditional framework of Galerkin projection-based model reduction but modifies the reduced order model so that the appropriate Hamiltonian structure is preserved. However, its computational complexity for online simulations is still high if the Hamiltonian involves non-polynomial nonlinearities. In this paper, we apply the discrete empirical interpolation method to improve the online efficiency of the structure-preserving reduced order simulations. Since the reduced basis truncation can degrade the Hamiltonian approximation, we propose to use the basis vectors obtained from shifted snapshots. 
A nonlinear wave equation is used as a test bed and the numerical results illustrate the efficacy of the proposed method. 
\end{abstract}



\newcommand{\lp}{\left(}
\newcommand{\rp}{\right)}
\newcommand{\lno}{\left\|}
\newcommand{\rno}{\right\|}
\newcommand{\id}{\text{ d}}

\newcommand{\ou}{\overline{u}}
\newcommand{\opsi}{\overline{\psi}}
\newcommand{\oq}{\overline{q}}
\newcommand{\oT}{\overline{T}}
\newcommand{\E}{\mathbbm{E}}
\newcommand{\orho}{\overline{\rho}}

\newcommand{\s}{\sigma}
\renewcommand{\k}{\kappa}
\newcommand{\p}{\partial}
\newcommand{\om}{\omega}
\newcommand{\Om}{\Omega}
\newcommand{\pOm}{\partial \Omega}
\newcommand{\e}{\epsilon}
\renewcommand{\a}{\alpha}
\renewcommand{\b}{\beta}
   \newcommand{\eps}{\varepsilon}
   \newcommand{\EX}{{\Bbb{E}}}
   \newcommand{\PX}{{\Bbb{P}}}

\newcommand{\nd}{{\nabla \cdot}}

\newcommand{\cF}{{\cal F}}
 
\newcommand{\cD}{{\cal D}}
\newcommand{\cO}{{\cal O}}

\newtheorem{remark}{Remark}[section]
\newtheorem{lemma}{Lemma}[section]
\newtheorem{theorem}{Theorem}[section]
\newtheorem{corollary}{Corollary}[section]
\newtheorem{proposition}{Proposition}[section]
\newtheorem{definition}{Definition}[section]
\newcommand{\kmodel}{k_{\mbox{Model}}}
\newcommand{\obu}{\overline{\bf u}}
\newcommand{\oobu}{\overline{\overline{\bf u}}}
\newcommand{\be}{{\bf e}}
\newcommand{\bk}{{\bf k}}
\newcommand{\bs}{{\bf s}}
\newcommand{\bu}{{\bf u}}
\newcommand{\bD}{{\bf D}}
\newcommand{\bN}{{\bf N}}
\newcommand{\bS}{{\bf S}}
\newcommand{\oou}{\overline{\overline{u}}}
\newcommand{\op}{\overline{p}}
\newcommand{\of}{\overline{f}}
\newcommand{\obf}{\overline{\bf f}}
\newcommand{\ow}{\overline{w}}
\newcommand{\ov}{\overline{v}}
\newcommand{\ophi}{\overline{\phi}}
\newcommand{\oS}{\overline{S}}
\newcommand{\obS}{\overline{\bf S}}
\newcommand{\bv}{{\bf v}}
\newcommand{\obv}{\overline{\bf v}}
\newcommand{\bc}{{\bf c}}
\newcommand{\by}{{\bf y}}
\newcommand{\bA}{{\bf A}}
\newcommand{\bB}{{\bf B}}
\newcommand{\bG}{{\bf G}}
\newcommand{\bI}{{\bf I}}
\newcommand{\bQ}{{\bf Q}}
\newcommand{\bY}{{\bf Y}}
\newcommand{\bw}{{\bf w}}
\newcommand{\bW}{{\bf W}}
\newcommand{\bU}{{\bf U}}
\newcommand{\obw}{\overline{\bf w}}
\newcommand{\bz}{{\bf z}}
\newcommand{\bZ}{{\bf Z}}
\newcommand{\obZ}{\overline{\bf Z}}
\newcommand{\bff}{{\bf f}}
\newcommand{\bee}{{\bf e}}
\newcommand{\bn}{{\bf n}}
\newcommand{\bx}{{\bf x}}
\newcommand{\bX}{{\bf X}}
\newcommand{\bH}{{\bf H}}
\newcommand{\bV}{{\bf V}}
\newcommand{\bL}{{\bf L}}
\newcommand{\bg}{{\bf g}}
\newcommand{\bj}{{\bf j}}
\newcommand{\br}{{\bf r}}
\newcommand{\grads}{\nabla^s}
\def\PP{{{\rm l}\kern - .15em {\rm P} }}
\def\PN2{{\PP_{N}-\PP_{N-2}}}
\newcommand{\erf}[1]{\mbox{erf}\left(#1\right)}
\newcommand{\D}{\nabla}
\newcommand{\I}{\mathbb{I}}
\newcommand{\N}{\mathbb{N}}
\newcommand{\R}{\mathbb{R}}
\newcommand{\Z}{\mathbb{Z}}
\newcommand{\mathR}{\R}
\newcommand{\mathN}{\N}
\newcommand{\mathZ}{\Z}
\newcommand{\mathI}{\mathbbm{I}}
\newcommand{\btau}{\boldsymbol{\tau}}
\newcommand{\bphi}{\boldsymbol{\phi}}
\newcommand{\bvarphi}{\boldsymbol{\varphi}}
\newcommand{\bpsi}{\boldsymbol{\psi}}
\newcommand{\bfeta}{\boldsymbol{\eta}}
\newcommand{\blambda}{\boldsymbol{\lambda}}
\newcommand{\bPhi}{\boldsymbol{\Phi}}
\newcommand{\bPsi}{\boldsymbol{\Psi}}
\newcommand{\obphi}{\overline {\boldsymbol{\phi}}}
\newcommand{\bomega}{\boldsymbol{\omega}}
\newcommand{\bsigma}{\boldsymbol{\sigma}}
\newcommand{\orhoprime}{\overline{\rho^{\prime}}}
\newcommand{\bus}{{\bf u}^*}
\newcommand{\By}{\mathcal B(\by)}
\newcommand{\eci}[1]{\mathcal E_{#1}}
\newcommand{\dpyi}[1]{\delta_{#1}^+(\by)}
\newcommand{\dmyi}[1]{\delta_{#1}^-(\by)}
\newcommand{\cA}{{\mathcal A(\by)}}
\newcommand{\dyi}[1]{\delta_{#1}(\by)}
\newcommand{\cG}{{\mathcal G(\bx,\by)}}
\newcommand{\cGi}[1]{{\mathcal G_{#1}(\bx,\by)}}
\newcommand{\pti}{\partial_i}
\newcommand{\ptii}[1]{\partial_{#1}}
\newcommand{\ba}{{\bf a}}

\newcommand{\rey}{\mbox{Re}}

\newcommand{\tnp}{t^{k+1}}
\newcommand{\bb}{{\bf b}}
\newcommand{\fnp}{f^{k+1}}
\newcommand{\prp}{P_R^{'}}
\newcommand{\pr}{P_R}
\newcommand{\rn}{r^k}
\newcommand{\ur}{\bu_r}
\newcommand{\urn}{\bu_r^k}
\newcommand{\utn}{\bu_t^k}
\newcommand{\urnp}{\bu_r^{k+1}}
\newcommand{\utnp}{\bu_t^{k+1}}
\newcommand{\un}{\bu^k}
\newcommand{\unp}{\bu^{k+1}}
\newcommand{\vr}{\bv_r}
\newcommand{\wrn}{\bw_r^k}
\newcommand{\etan}{\eta^k}
\newcommand{\etanp}{\eta^{k+1}}
\newcommand{\prn}{\phi_r^k}
\newcommand{\prnp}{\phi_r^{k+1}}
\newcommand{\prz}{\phi_r^{0}}
\newcommand{\prm}{\phi_{r}^{M}}
\newcommand{\sNn}{\sum\limits_{k=0}^{M-1}}
\newcommand{\sN}{\sum\limits_{k=0}^{M}}
\newcommand{\asNn}{\frac{1}{M-1}\sNn}
\newcommand{\asN}{\frac{1}{M}\sN}
\newcommand{\uph}{\upsilon^h}
\newcommand{\half}{\frac{1}{2}}
\section{Introduction}
\noindent \indent There are broad applications of Hamiltonian systems in engineering and scientific research.  For real-world problems, using numerical methods to simulate such systems often requires long-time integrations of large-scale discrete systems. Because the numerical error accumulates over time, preserving intrinsic properties of the systems has been a key criterion for developing stable numerical schemes.
So far, geometric integrators or structure-preserving algorithms have been introduced to exactly preserve structural properties of Hamiltonian systems. For instance, symplectic algorithms proposed in \cite{Hairer02,Feng03} has achieved a remarkable success in dealing with Hamiltonian ordinary differential equations (ODEs). They have been extended to Hamiltonian partial differential equations (PDEs) and preserve the multi-symplectic conservation law \cite{Bridges06,Wang2013}.
In recent years, there has been an increasing emphasis on constructing numerical methods to
preserve certain invariant quantities such as the total energy of dynamical systems.
Several discrete gradient methods have been proposed in the literature 
\cite{Gonzalez96,McLachlan99,Quispel08,Celledoni2012,Gong2014}.


To accelerate the long-time, large-scale numerical simulations of the Hamiltonian systems, reduced-order modeling has been considered. 
One such model reduction technique is the proper orthogonal decomposition (POD) method, which has been successfully applied to many time-dependent, nonlinear PDEs (\cite{bui2007goal,carlberg2011low,chaturantabut2010nonlinear,daescu2008dual,HLB96,iollo2000stability,KV01,sirisup2004spectral,lassila2014model}).
The POD method extracts orthogonal basis vectors from snapshot data, and use them to span the trial space. 
Based on different choices of test space, either Galerkin or Petrov-Galerkin projection can be used to build a reduced order model (ROM). Such a ROM is low-dimensional, but it does not possess the Hamiltonian structure. 
This issue has been recently recognized and structure-preserving ROMs (SP-ROMs) have been developed to resolve it. 
For instance, structure-preserving Petrov-Galerkin reduced models were introduced in \cite{beattie2011structure,chaturantabut2016structure} for port-Hamiltonian systems, in which the POD-based and $\mathcal{H}_2$-based projection subspaces were considered. 
A proper symplectic decomposition approach using the symplectic Galerkin projection was proposed in \cite{peng2016symplectic} for Hamiltonian PDEs with a symplectic structure. 
A structure-preserving POD-ROM was introduced in \cite{gong2017structure}, where the Galerkin projection-based ROM was modified so that appropriate Hamiltonian structure can be well kept. The approach has been applied to nonlinear Schr\"{o}dinger equation and shallow water equations in \cite{karasozen2018energy,karasozen2021structure,sockwell2019mass}. 
Meanwhile, when non-polynomial nonlinearities appear in SP-ROMs, their computational costs still depend on the number of spatial degrees of freedom in the full-order model. 
Thus, a hyper-reduction method such as the discrete empirical interpolation method (DEIM) \cite{chaturantabut2012state} has been employed as a remedy. 
It was combined with the Kronecker product in \cite{miyatake2019structure} to effectively reduce the complexity for evaluating the variable skew-symmetric coefficient matrix. In \cite{chaturantabut2016structure,karasozen2018energy}, two different ways for applying DEIM to the gradient of Hamiltonian functions were proposed, respectively. 
It is worth mentioning that other types of model reduction techniques such as reduced basis method have been introduced for Hamiltonian systems as well \cite{afkham2017structure,afkham2019structure,hesthaven2020rank,pagliantini2020dynamical}; and ROMs have been investigated for preserving different geometric properties such as the Lagrangian structures in \cite{carlberg2015preserving,farhat2015structure}.

In this paper, we consider the framework of the SP-ROM developed in \cite{gong2017structure} and extend it to the case in which the gradient of Hamiltonian involves non-polynomial nonlinearities. A structure-preserving POD-DEIM ROM is then proposed that possesses the appropriate Hamiltonian structure while reducing the online simulation cost. Due to the basis truncation, there may exist discrepancies in the Hamiltonian between the new ROM and the full-order model.
Therefore, we propose to use the POD basis and DEIM basis from shifted snapshots to improve the Hamiltonian approximation.

The rest of this paper is organized as follows.
In Section \ref{sec: alg}, we introduce the SP-ROMs for Hamiltonian systems including the SP-POD and SP-DEIM models;
The SP-ROMs are numerically investigated in Section \ref{sec: num};
A few concluding remarks are drawn in the last section.

\section{Structure-Preserving Galerkin ROMs}\label{sec: alg}

\noindent \indent Consider a general Hamiltonian PDE system
\begin{equation}
\dot{\bu} = \mathcal{D}\, \frac{\delta \mathcal{H}}{\delta \bu},
\label{eq: ham_pde}
\end{equation}
where $\mathcal{D}$ is a differential operator and $\mathcal{H}(\bu)$ is the Hamiltonian, which often corresponds to the total energy of the system.
When $\mathcal{D}$ is a skew-adjoint operator with respect to the $L_2$ inner product, the PDE keeps $\mathcal{H}(\bu)$ invariant, which can be easily checked using the fact of  
$\int_{\Omega} \frac{\delta \mathcal{H}}{\delta \bu} \mathcal{D} \frac{\delta \mathcal{H}}{\delta \bu} dx = 0$. 
Note that when $\mathcal{D}$ is a constant negative semi-definite (or definite) operator with respect to the $L_2$ inner product, the system becomes dissipative and $\mathcal{H}(\bu)$, referred to as the Lyapunov function, would be non-increasing. In this paper, we are concerned with the former case, although the extension to the latter is natural.

Since \eqref{eq: ham_pde} preserves the Hamiltonian when $\mathcal{D}$ is skew-adjoint, numerical methods for the PDE are expected to keep the same property at a discrete level. This has become a basic rule of thumb when designing a robust numerical scheme for the Hamiltonian PDE, especially for the purpose of long term simulations. Therefore, methods such as geometric integrators \cite{Hairer02,Feng03,Bridges06,Wang2013} or structure-preserving algorithms \cite{Quispel08,Hairer10,Cohen11,Celledoni2012,Dahlby2011,Gong2014} have been developed to preserve such structural properties. 
On the other hand, for large-scale numerical simulations, in order to reduce the computational complexity, structure-preserving reduced order modeling has been introduced. 
Next, we first briefly review the structure-preserving approach proposed in \cite{gong2017structure}, and then extend it to general nonlinear cases using the DEIM.

\subsection{The structure-preserving POD-ROM}\label{sec: sp-pod}
\noindent \indent Assume a suitable spatial discretization of (\ref{eq: ham_pde}) yields a finite $n$-dimensional Hamiltonian ODE system of the form  (e.g. see \cite{Celledoni2012}):
\begin{equation}
\dot{\bu} = \bD\, \nabla_{\bu} H(\bu)
\label{eq: ham_fom}
\end{equation}
with the initial condition $\bu(t_0)= \bu_0$ and appropriate boundary conditions, where the coefficient matrix $\bD\in \mathbb{R}^{n\times n}$ on the right-hand-side is skew symmetric. 
In our numerical experiments, we use the finite difference method for spatial discretization, in which $n$ frequently equals to the number of grid points; and $H: \mathbb{R}^n \rightarrow \mathbb{R}$ such that $H\Delta x$ provides a consistent approximation to the Hamiltonian $\mathcal{H}(\bu)$ associated to \eqref{eq: ham_pde}. 
Note that we slightly abuse the notation here: $\bu$ is used again to denote the discrete state variable in $\mathbb{R}^n$.

Given the reduced basis matrix $\bPhi\in \mathbb{R}^{n\times r}$, the state variable $\bu$ can be approximated by $\bu_r(t) = \bPhi \ba(t)$, where $\ba(t)$ is the unknown, $r$-dimensional coefficient vector. The reduced basis vectors can be determined at an offline stage by different model reduction techniques such as the reduced basis methods, proper orthogonal decomposition (POD), dynamical mode decomposition, etc. Here, we choose the POD because it is relatively easy to use and has been applied to many practical engineering problems. 
A standard Galerkin projection-based reduced order model (G-ROM) reads: 
\begin{equation}
\dot{\ba} = \bPhi^\top\bD\, \nabla_{\bu} H(\bPhi \ba).
\label{eq: ham_rom1}
\end{equation}
It provides a low-dimensional surrogate to \eqref{eq: ham_fom} since one only solves $\ba(t)$ at the online stage and $\bu_r$ can be computed if needed once $\ba(t)$ is determined. However, this model does not provide a constant Hamiltonian approximation, which can be seen by checking the time derivate of $H_r(t) \equiv H(\bPhi \ba(t))$ as follows. 
\begin{eqnarray*}
\frac{d}{dt}H(\bPhi \ba)&=&[\nabla_\ba H(\bPhi \ba)]^\top \dot{\ba} \nonumber \\
				  &=&[\bPhi^\top \nabla_{\bu} H(\bPhi \ba)]^\top \bPhi^\top\bD\, \nabla_{\bu} H(\bPhi \ba) \nonumber \\
				  &=& \nabla_{\bu} H(\bPhi \ba)^\top \bPhi \bPhi^\top\bD\, \nabla_{\bu} H(\bPhi \ba),
\end{eqnarray*}
where we use the fact that
$
\nabla_{\ba} H(\bPhi \ba) = \bPhi^\top \nabla_{\bu} H(\bPhi \ba)
$.
Since $\bPhi$ is composed of the first $r$ left singular vectors of the snapshot matrix, $\bPhi\bPhi^\top\bD$ is not skew symmetric as $\bD$, hence, the time derivative of $H(\bPhi \ba)$ is not ensured to be zero.

The structure-preserving POD (SP-POD) reduced order model developed in \cite{gong2017structure} modifies the above approach such that the Hamiltonian is well preserved, which has the following form: 
\begin{equation}
\dot{{\ba}} = \bD_r\, \nabla_{\ba} H(\bPhi {\ba}),
\label{eq: ham_rom2}
\end{equation}
where 
$
\bD_r = \bPhi^\top \bD \bPhi 
$. 
For this reduced-order dynamical system, the time derivative of $H_r(t)$ is
\begin{equation}
	\begin{aligned}
	\frac{d}{dt}H(\bPhi  {\ba})&=[\nabla_{ {\ba}} H(\bPhi {\ba})]^\top \dot{ {\ba}} \nonumber \\
				  &=[\nabla_{ {\ba}} H(\bPhi  {\ba})]^\top \bD_r\, \nabla_{\ba} H(\bPhi  {\ba}) \nonumber \\
				  &= 0,\quad \text{ as $\bD_r$ is skew-symmetric.}
	\end{aligned}
\label{eq: dHdt}	
\end{equation}
Hence, this model possesses an invariant $H_r(t)$. With a structure-preserving time stepping, the discrete Hamiltonian $H_r \Delta x$ can be well preserved during long-term simulations.

\subsection{The structure-preserving DEIM-ROM}\label{sec: sp_deim1}

Although \eqref{eq: ham_rom2} has a low dimension, its computational cost could still depend on the number of spatial degrees of freedom if $\bD= \bD(\bu)$ depends nonlinearly on $\bu$ or $\nabla_{\bu} H(\bu)$ includes non-polynomial nonlinearities. In either case, online simulations involve certain calculations to be performed at all the grid points. 
To overcome this issue, hyper-reduction methods such as the DEIM have been applied. 
For instance, DEIM is combined with tensor product and vectorization in \cite{miyatake2019structure} to effectively reduce the complexity for evaluating $\bD_r(\bu_r)$. 
We shall focus on the case in which $\nabla_{\bu} H(\bu)$ has non-polynomial nonlinearities in this work. 
For completeness of presentation, we next first present the DEIM algorithm, review existing DEIM Hamiltonian approximations, and then introduce a new approach.    

In general, the DEIM employs the following ansatz on a nonlinear function $f(\bu(t))$:
\begin{equation}
f(\bu(t)) = \sum\limits_{j=1}^{s} \bpsi_j c_j(t),
\end{equation}
where $\bpsi_j$ is the $j$-th DEIM basis vectors generated from the nonlinear snapshots
$$[f(\bu(t_1)), f(\bu(t_2)), \ldots, f(\bu(t_m))].$$ 
The DEIM, as shown in Algorithm \ref{alg: DEIM}, selects a set of interpolation points $\wp := [\wp_1, \ldots, \wp_s]^{\intercal}$ in a greedy manner, in which $e_{\wp_i}$ be the $\wp_i$-th column of the identity matrix.

\begin{algorithm}\caption{DEIM}\label{alg: DEIM}
\SetKwInOut{Input}{input}\SetKwInOut{Output}{output}
\vspace{.3cm}
\Input{$\{\bpsi_{\ell}\}_{\ell=1}^{s} \subset \mbr^{s}$ linear independent}
\Output{$\wp = [\wp_1, \ldots, \wp_s]^{\intercal} \in \mbr^s$}
$[|\rho|,\, \wp_1] = \max\{|\bpsi_1|\}$\;
$\bPsi = [\bpsi_1], \bP = [{\bf e}_{\wp_1}], \wp = [\wp_1]$\;
\For{$\ell = 2$ \KwTo $s$}{
Solve $(\bP^{\intercal} \bPsi){\bf c} = \bP^{\intercal} \bpsi_{\ell}$ for $\bf c$ \;
${\bf r}=\bpsi_{\ell}- \bPsi {\bf c}$\;
$\left[ |\rho|, \wp_{\ell} \right] = \max\{|{\bf r}|\}$\;
$\bPsi \leftarrow [\bPsi \quad \bpsi_{\ell}], \bP\leftarrow [\bP\quad {\bf e}_{\wp_{\ell}}], \wp \leftarrow
\left[\begin{array}{c} \wp \\ \wp_{\ell}\end{array}\right]$\;
}
\end{algorithm}

The DEIM approximation of the nonlinear function is given by
\begin{equation}
f(\bu) \approx {\bf \Psi}(\bP^\intercal {\bf \Psi})^{-1} \bP^\intercal {f}(\bu),
\label{eq:deim}
\end{equation}
where $\bP = [e_{\wp_1}, \ldots, e_{\wp_s}]\in \mathbb{R}^{n\times s}$ and ${\bf \Psi} \in \mathbb{R}^{n\times s}$ is the DEIM basis matrix. Because the DEIM approximation only evaluate $f$ at $s$ points, it could greatly reduce the online computation. 
For a detailed description of the method, the read is referred to \cite{chaturantabut2010nonlinear}. 

When the discrete Hamiltonian contains non-polynomial nonlinear functions of $\bu$, the gradient $\nabla_{\bu} H(\bu)$ usually includes non-polynomial nonlinearities as well. It makes the online computational complexity still depend on $n$, the number of degrees of freedom of the full-order model. To overcome this issue, a DEIM hamiltonian was proposed in \cite{chaturantabut2016structure}. 
It splits $H(\bu)$ into two parts: 
\begin{equation}
H(\bu) = \frac{1}{2} \bu^\intercal \bQ \bu + h(\bu).
\label{eq: ham}
\end{equation}
The first part is quadratic in $\bu$, where $\bQ\in \mathbb{R}^{n\times n}$ is constant and positive definite, which contributes the linear component of $\nabla_{\bu} H(\bu)$; the second one is the remainder, which is nonlinear and yields the nonlinearity of $\nabla_{\bu} H(\bu)$.  
It is natural to approximate $\nabla_{\bu} h(\bu)$ by a DEIM interpolation: 
$\nabla_{\bu} h(\bu(t)) \approx \bPsi (\bP^\intercal \bPsi)^{-1} \bP^\intercal\nabla_{\bu} h(\bu(t)) $. 
Denote the DEIM projection by $\mathbb{P} = \bPsi (\bP^\intercal \bPsi)^{-1} \bP^\intercal$, the DEIM Hamiltonian is defined in \cite{chaturantabut2016structure} by   
$$\widehat{H}(\bu) = \frac{1}{2} \bu^\intercal \bQ \bu + h(\mathbb{P}^\intercal \bu),
$$
whose gradient is 
$$\nabla_{\bu}\widehat{H}(\bu) = \bQ \bu + \mathbb{P} \nabla_{\bu} h(\mathbb{P}^\intercal \bu).
$$
Then it can be shown that $\widehat{H}_r(t)\equiv \widehat{H}(\bPhi \ba(t))$ is invariant. 
However, $\bu$ is approximated by $\mathbb{P}^\intercal \bu$ in this approach but $\mathbb{P}$ is derived from the DEIM on $\nabla_{\bu} h(\bu)$, which could introduce errors.
In \cite{karasozen2018energy}, a nonlinear function in $\nabla_{\bu} H(\bu)$ denoted by $\nabla_{\bu} h(\bu)$ is directly approximated by its DEIM interpolation $\mathbb{P} \nabla_{\bu} h(\bu)$. This way could avoid approximating $\bu$ by $\mathbb{P}^\intercal \bu$, but there does not exist an explicit formula of the discrete Hamiltonian.  
As a result, this approach cannot preserve the discrete energy exactly.

%

Next, we propose a new way to deal with the nonlinearity in the gradient of $H(\bu)$ so that  
the computational complexity of nonlinear function evaluations is independent of $n$ while keeping the Hamiltonian invariant at the discrete level. 
The key idea is to apply DEIM to the Hamiltonian function, not to its gradient. 
To this end, we first recognize a discrete Hamiltonian function can usually be written as follows:  
$$
H(\bu) = \frac{1}{2} \bu^\intercal \bQ \bu + \bc^\intercal \bG (\bu), 
$$
where $\bG$ is a nonlinear vector-valued function of $\bu$, and $\bc^\intercal \bG (\bu)$ represents the non-polynomial nonlinearity. 
After collecting snapshots of $\bG(\bu)$, $[\bG(\bu(t_0), \ldots, \bG(\bu(t_m)]$, we generate the following DEIM interpolation for $\bG(\bu)$ by Algorithm \ref{alg: DEIM}:
\begin{equation*}
\bG(\bu) \approx \mathbb{P} \bG(\bu) = {\bPsi}(\bP^\intercal {\bPsi})^{-1} \bP^\intercal \bG(\bu).
\label{eq:deim_g}
\end{equation*}
The reduced Hamiltonian is then defined by  
\begin{equation}
H_r(\bPhi\ba) = \frac{1}{2} \ba^\intercal \bQ_r \ba + \bc^\intercal \mathbb{P} \bG (\bPhi \ba), 
\label{eq:deim_ham}
\end{equation}
with $\bQ_r = \bPhi^\intercal \bQ \bPhi$ and the associated gradient can be expressed as 
\begin{equation}
\nabla_{\ba}H_r(\ba) = \bQ_r \ba + \bPhi^\intercal \bJ_{\bG} (\bPhi \ba) \mathbb{P}^\intercal \bc, 
\label{eq:deim_gham}
\end{equation}
where $\bJ_{\bG}(\cdot)$ is the Jacobian matrix of $\bG(\bu)$. 
Plugging \eqref{eq:deim_gham} into \eqref{eq: ham_rom2}, we have a structure-preserving DEIM (SP-DEIM) model: 
\begin{equation}
\begin{aligned}
\dot{{\ba}} &= \bD_r\, \nabla_{\ba} H_r(\bPhi {\ba}) \\
			&= \bD_r\, \left[\bQ_r \ba + \bPhi^\intercal \bJ_{\bG} (\bPhi \ba) \mathbb{P}^\intercal \bc\right]. 
\label{eq: ham_rom3}
\end{aligned}
\end{equation}
Using the same argument as \eqref{eq: dHdt}, we can easily show the discrete Hamiltonian $H_r$ is a constant. 

\begin{remark} The philosophy of the proposed approach to treat the nonlinear Hamiltonian is similar to that of the way to treat nonlinear Lagrangian in \cite{carlberg2015preserving}: the reduction and hyper-reduction are first performed to approximate the key quantity, Lagrangian in \cite{carlberg2015preserving} or Hamiltonian in this work, and the SP-ROMs are then derived based on the reduced quantity. But \cite{carlberg2015preserving} is concerned with second-order dynamics, while we consider general energy-preserving systems. For classical mechanical systems, our approach could lead to SP-ROMs equivalent to reduced order Euler-Lagrangian equations. For instance, for the nonlinear wave equation tested in Section \ref{sec: SP-DEIM}, if the basis of $u$ is used for $v$ as well, the SP-DEIM has an equivalent second-order reduced-order dynamics that can be derived from \cite{carlberg2015preserving}. 
\end{remark} 
\subsection{Improvement of the Hamiltonian approximation}\label{sec: sp_deim2}
A small amount of POD basis functions captures most of the snapshot information, but the basis truncation would result in a loss of information. It could further result in a discrepancy in the Hamiltonian function approximation between the ROM and the FOM.
Since the proposed structure-preserving ROMs are able to keep a constant Hamiltonian, it is enough to correct the initial value of the discrete Hamiltonian to improve the approximation during the reduced-order simulations.
Therefore, we follow one approach developed in \cite{gong2017structure} that uses the POD basis generated from shifted snapshots. 

Instead of the snapshots of state variables, $\left[\bu(t_1), \ldots, \bu(t_M)\right]$, we consider a shifted snapshot set
\begin{equation}
\left[\bu(t_1)-\bu_0, \ldots, \bu(t_M)-\bu_0\right]
\label{eq: snap_shift}
\end{equation}
and obtain the POD basis $\bPhi$. The associated reduced approximation becomes
\begin{equation}
\bu_r(t)= \bPhi {\ba}(t) + \bu_0. 
\label{eq: pod_shift}
\end{equation}
Accordingly, we generate DEIM basis $\bPsi$ from a shifted nonlinear snapshot set 
\begin{equation}
\left[\bG(\bu(t_1))-\bG(\bu_0), \ldots, \bG(\bu(t_M))-\bG(\bu_0) \right].
\label{eq: nonlsnap_shift}
\end{equation}
The related DEIM approximation becomes
$
\mathbb{P} (\bG(\bu)-\bG(\bu_0)) + \bG(\bu_0)
$
with $\mathbb{P}$ the DEIM projection. The reduced Hamiltonian is 
\begin{equation}
H_r(t) = \frac{1}{2} (\bu_0^\intercal + \ba^\intercal \bPhi^\intercal) \bQ (\bPhi \ba+\bu_0) 
+ \bc^\intercal [\mathbb{P}\bG (\bPhi \ba+\bu_0) + (\mathbb{I}-\mathbb{P}) \bG(\bu_0)]. 
\label{eq:deim_ham_new}
\end{equation}
Substituting \eqref{eq:deim_ham_new} into \eqref{eq: ham_rom2}, we have
\begin{equation}
\begin{aligned}
\dot{{\ba}} &= \bD_r\, \nabla_{\ba} H_r(\bPhi \ba + \bu_0) \\
			&= \bD_r\, \left[\bQ_r \ba + \bPhi^\intercal \bQ \bu_0 + \bPhi^\intercal \bJ_{\bG} (\bPhi \ba + \bu_0) \mathbb{P}^\intercal \bc\right], 
\label{eq: ham_rom3}
\end{aligned}
\end{equation}
together with the initial condition ${\ba}(t_0) = \bf{0}$.
Since $\bD_r$ is skew symmetric, the model preserves the structural property and the reduced Hamiltonian would keep the same value as $H(\bu_0)$ during the reduced order simulations. 

\section{Numerical experiments\label{sec: num}}
\noindent \indent In this section, we use a nonlinear wave equation as an example of Hamiltonian PDEs to investigate the numerical performance of the structure-preserving ROMs. 
Consider a one-dimensional case with a constant moving speed $c$ and a nonlinear forcing term $g(u)$,
\begin{equation*}
u_{tt}= c^2 u_{xx} - g(u), \quad 0\leq x\leq l.
\end{equation*}
The equation can be written in the Hamiltonian formulation
\begin{equation}
\left[
\begin{array}{c}
\dot{u}\\
\dot{v}
\end{array}
\right] =
\left[
\begin{array}{cc}
0 & 1 \\
-1 & 0
\end{array}
\right]
\left[
\begin{array}{c}
\frac{\delta \mathcal{H}}{\delta u}\\
\frac{\delta \mathcal{H}}{\delta v}
\end{array}
\right],
\label{eq:lin_wave}
\end{equation}
which has a symplectic structure. The PDE system has a constant Hamiltonian that is also the total energy:  
\begin{equation*}
\mathcal{H}(t) = \int_0^l \left[\frac{1}{2}v^2+\frac{c^2}{2} u_x^2 + G(u)\right]\, dx,
\end{equation*}
where ${G}'(u) = g(u)$,
$\frac{\delta \mathcal{H}}{\delta u} = -c^2 u_{xx}+g(u)$ and $\frac{\delta \mathcal{H}}{\delta v} = v$.
After a spatial discretization using finite difference method with $n$ uniformly distributed grid points and the mesh size $\Delta x$, the semi-discrete system of \eqref{eq:lin_wave} reads: 
\begin{equation}
\left[
\begin{array}{c}
\dot{\bu}\\
\dot{\bv}
\end{array}
\right] =
\left[
\begin{array}{cc}
0 & \bI_n \\
-\bI_n & 0
\end{array}
\right]
\left[
\begin{array}{c}
-\bA \bu+\bg(\bu)\\
\bv
\end{array}
\right],
\label{eq: wavef}
\end{equation}
where $\bA$ is a discrete, scaled, one-dimensional second order differential operator.
The discrete Hamiltonian is $H\Delta x$ with 
$$
H(t) = \frac{1}{2} \bv^\intercal \bv - \frac{1}{2} \bu^\intercal \bA \bu + \bc^\intercal \bG(\bu),
$$   
where $\bc$ is the all-ones vector and $\bc^\intercal {\bG}({\bu})$ corresponds to the discretization of the integration of ${G}(u)$. 

The linear problem has been tested in \cite{peng2016symplectic,gong2017structure}. Here we focus on the nonlinear case $g(u) = \sin(u)$ and ${G}({u})= 1-\cos(u)$, then $[\bG(\bu)]_i= 1-\cos(u_i)$ and $[\bg(\bu)]_i = \sin(u_i)$ for $i=0, \ldots, n-1$.
In particular, $c= 0.1$, $x\in [0, 1]$, $t\in [0, 50]$, and periodic boundary conditions are selected in our setting.
The initial condition satisfies $u(0)=f(s(x))$ and $\dot{u}(0)=0$, where $f(s)$ is a cubic spline function defined by
\begin{equation*}
f(s) =
\left\{
\begin{array}{ll}
1-\frac{3}{2}s^2+\frac{3}{4}s^3 & \text{\quad if  \,\,} 0\leq s\leq 1, \\
\frac{1}{4}(2-s)^3 		       & \text{\quad if  \,\,} 1< s\leq 2, \\
0					       & \text{\quad if  \,\,} s> 2,
\end{array}
\right.
\end{equation*}
and $s(x)= 10|x-\frac{1}{2}|$.

In the full-order simulation, the spatial domain is partitioned into $n=500$ equal subintervals, thus the mesh size $\Delta x= 2\times 10^{-3}$.
We use the symplectic midpoint method for the time integration with the step size $\Delta t= 0.01$.
A three-point stencil finite difference method is taken for the spatial discretization of the 1D second order differential operator.
The full order solution at time $t_{k+1}$, $\bu_h^{k+1}$ and $\bv_h^{k+1}$, satisfies
\begin{equation*}
\left[
\begin{array}{c}
\frac{\bu_h^{k+1}-\bu_h^k}{\Delta t}\\
\frac{\bv_h^{k+1}-\bv_h^k}{\Delta t}
\end{array}
\right] =
\left[
\begin{array}{cc}
0 &\bI_n \\
-\bI_n & 0
\end{array}
\right]
\left[
\begin{array}{c}
-\bA\frac{\bu_h^{k+1}+\bu_h^k}{2} + \bg\left(\frac{\bu_h^{k+1}+\bu_h^k}{2}\right)\\
\frac{\bv_h^{k+1}+\bv_h^k}{2}
\end{array}
\right],
\end{equation*}
where
\begin{equation*}
\bA= \frac{c^2}{\Delta x^2}
\left( \begin{array}{cccccc}
-2 & 1 & 0  & 0 & \cdots & 1 \\
1 & -2 & 1 & 0 & \cdots & 0 \\
   &     & \ddots & \ddots & \ddots &  \\
 0   &    \cdots&0          &  1 & -2  & 1 \\
 1 &    \cdots      &0     &  0  &  1 & -2
\end{array} \right),
\end{equation*}
the initial data $\bu_h^0$ has the $i$-th component equals $h(s(x_i))$ for $1\leq i\leq n$ and $\bv_h^0= {\bf 0}$.
The nonlinear system is solved by the Picard iteration.  

\begin{figure}[htb]
\centering
\begin{minipage}[ht]{0.31\linewidth}
\includegraphics[width=1\textwidth]{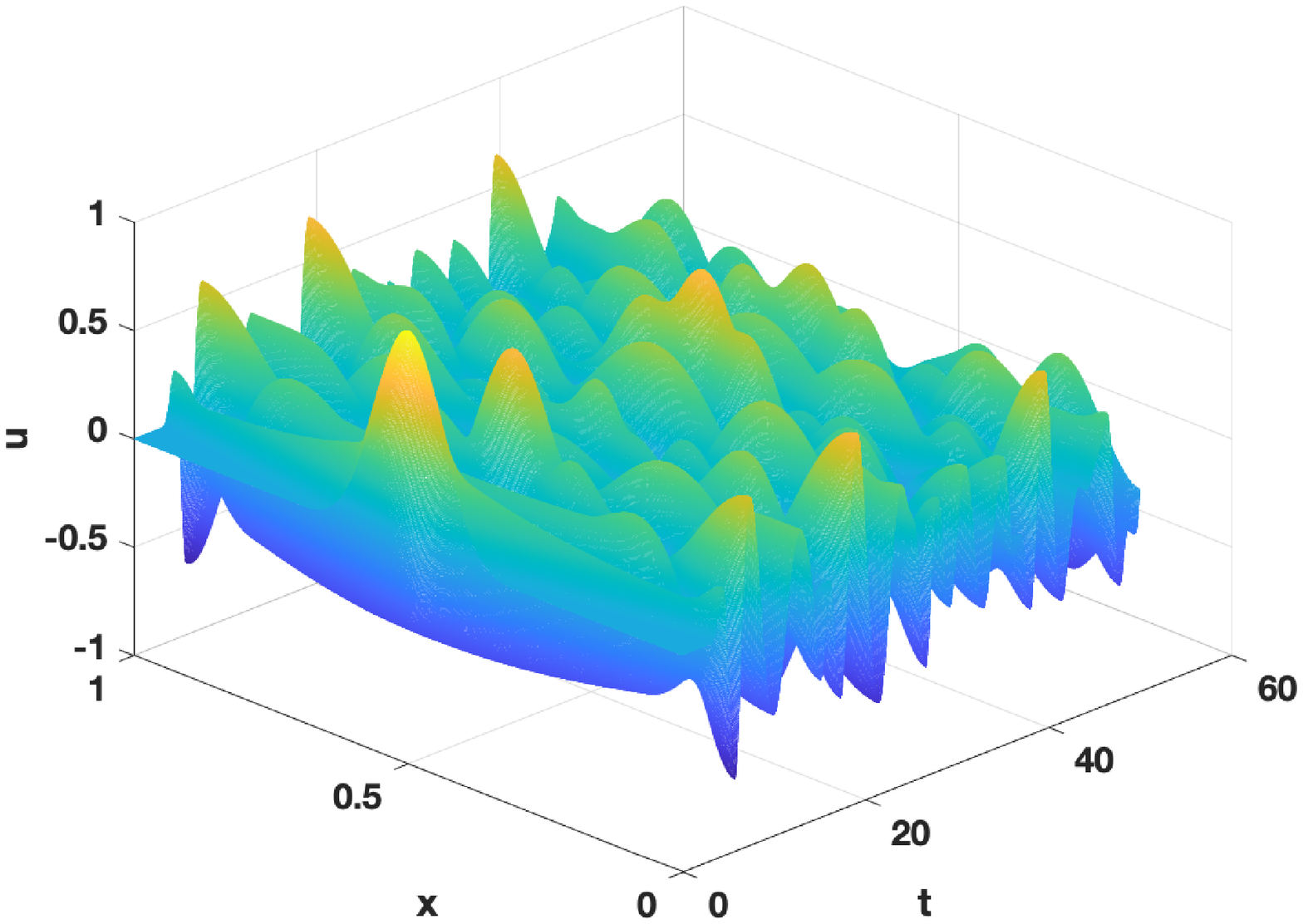}
\end{minipage}
\begin{minipage}[ht]{0.31\linewidth}
\includegraphics[width=1\textwidth]{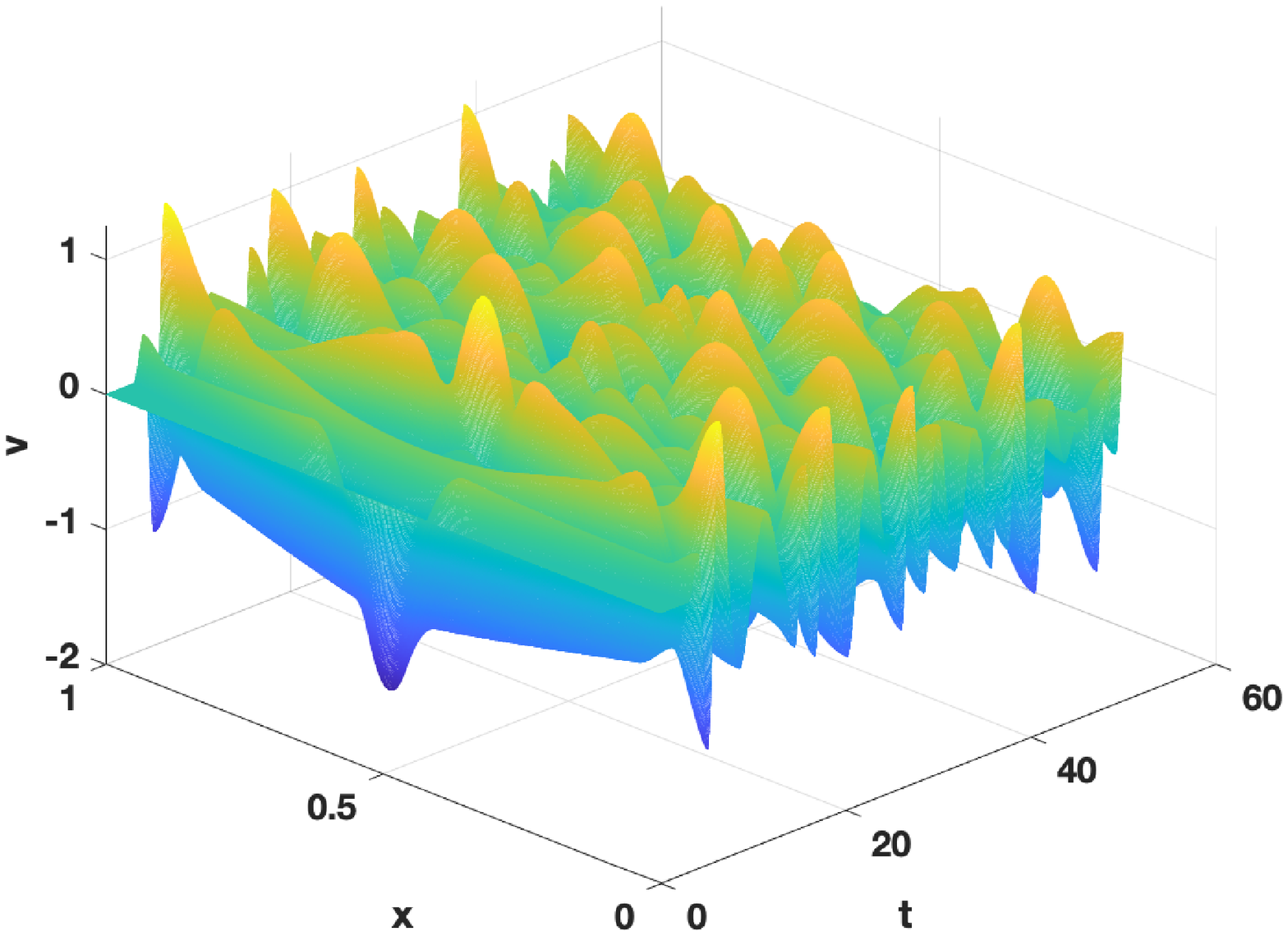}
\end{minipage}
\begin{minipage}[ht]{0.36\linewidth}
\vspace{1cm}
\includegraphics[width=1\textwidth]{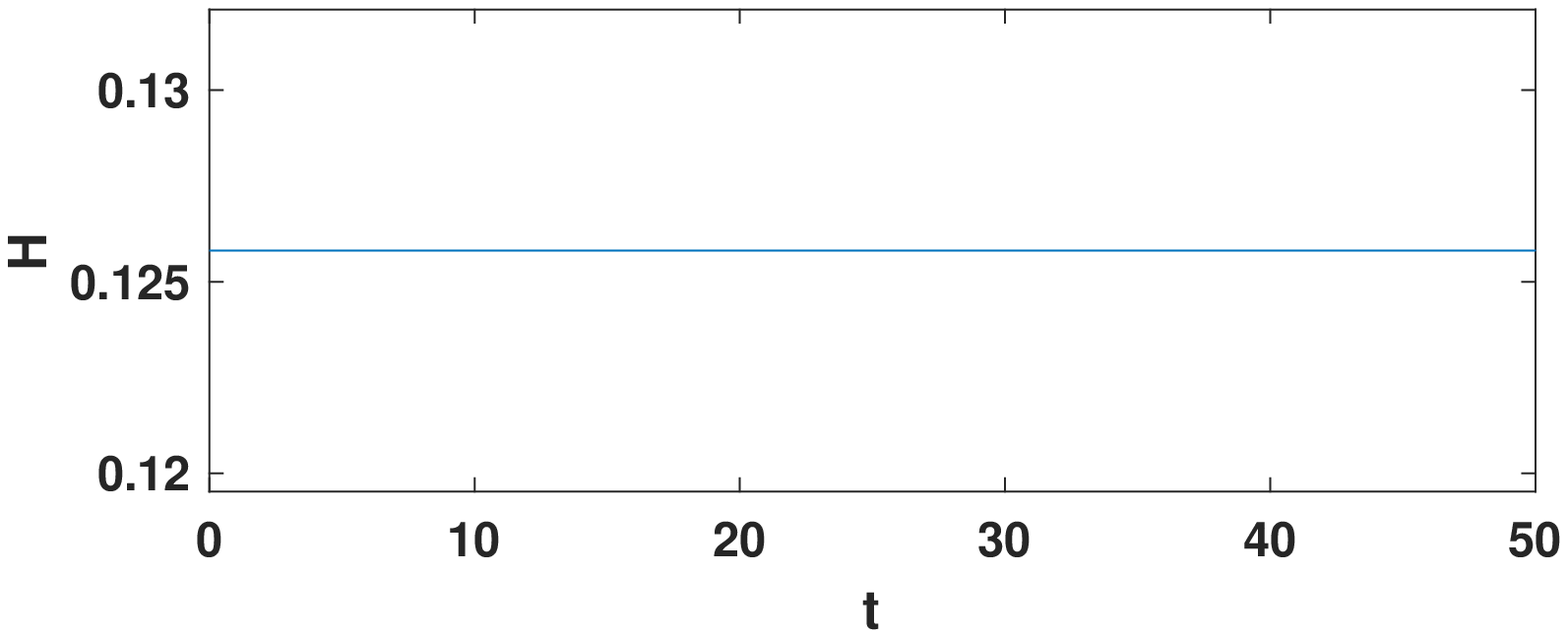}
\end{minipage}
\caption{
Full-order simulation results: time evolution of $\bu_h$ (left), $\bv_h$ (middle), and discrete Hamiltonian $H\Delta x$ (right).
}\label{Fig: lin_wave_full}
\end{figure}

It takes 52.8~seconds to finish the full-order simulation. The time evolution of numerical solutions $\bu_h$ and $\bv_h$ and the discrete energy $H\Delta x$ are plotted in Figure \ref{Fig: lin_wave_full}.
It is seen that the discrete energy is around $1.258\times 10^{-1}$.  
Since the exact solution is unknown, in what follows, the full-order simulation results will be served as benchmark solution for ROMs.

Next, we investigate the numerical performance of three types of ROMs: i) the G-POD model; ii) the SP-POD ROMs; iii) the SP-DEIM ROMs.
The criteria we shall use include: the maximum approximation error over the entire spatial-temporal domain,
$$\mathcal{E}_\infty= \max\limits_{k\geq 0} \max\limits_{0\leq i\leq n}  \sqrt{\left[(\bu_h^k)_i-(\bu_r^k)_i\right]^2+\left[(\bv_h^k)_i-(\bv_r^k)_i\right]^2},$$
the value of reduced-order Hamiltonian $H_r\Delta x$, 
and the CPU time $t_{cpu}$ for online simulations.

\subsection{Standard G-POD ROM.}
\noindent \indent Snapshots are collected from the full-order simulation every $50$ time steps.
Using the singular value decomposition, we find the $r$-dimensional POD basis $\bPhi_u$ and $\bPhi_v$ for $u$ and $v$, respectively.
The standard G-POD ROM is generated by substituting the POD approximation $\bu_r(t) = \bPhi_u \ba(t)$ and $\bv_r(t) = \bPhi_v \bb(t)$ into \eqref{eq: wavef} and applying the Galerkin projection, which can be written as: 
\begin{equation*}
\left[
\begin{array}{c}
\dot{\ba}\\
\dot{\bb}
\end{array}
\right] =
\left[
\begin{array}{c}
\bPhi_u^\top\bPhi_v \bb\\
 \bPhi_v^\top\bA \bPhi_u \ba - \bPhi_v^\intercal \bg(\bPhi_u \ba)
\end{array}
\right].
\label{eq:lin_wave_rom}
\end{equation*}
Using the same time integration method as the full order model, we have the POD basis coefficient at $t_{k+1}$, $\ba^{k+1}$ and $\bb^{k+1}$, satisfying
\begin{equation*}
\left[
\begin{array}{c}
\frac{\ba^{k+1}-\ba^{k}}{\Delta t}\\
\frac{\bb^{k+1}-\bb^{k}}{\Delta t}
\end{array}
\right] =
\left[
\begin{array}{c}
\bPhi_u^\top\bPhi_v \frac{\bb^{k+1}+\bb^{k}}{2}\\
\bPhi_v^\top\bA \bPhi_u \frac{\ba^{k+1}+\ba^{k}}{2}- \bPhi_v^\intercal \bg(\bPhi_u \frac{\ba^{k+1}+\ba^{k}}{2})
\end{array}
\right]
\label{eq:lin_wave_rom_dis}
\end{equation*}
with the initial condition $\ba^0 = \bPhi_u^\top \bu_0$ and $\bb^0 = \bPhi_v^\top \bv_0$.

\begin{figure}[htb]
\centering
\begin{minipage}[ht]{0.45\linewidth}
\includegraphics[width=1\textwidth]{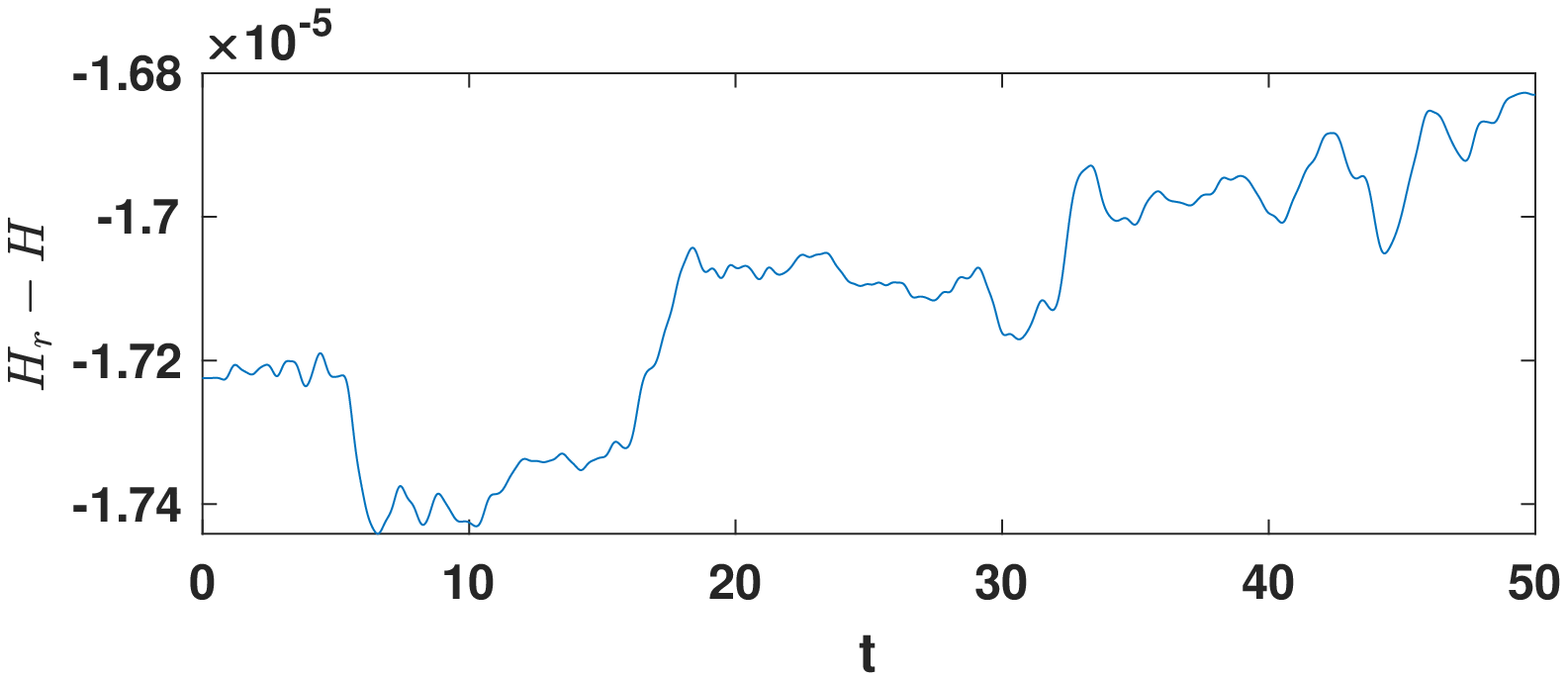}
\end{minipage}
\hspace{.2cm}
\begin{minipage}[ht]{0.45\linewidth}
\includegraphics[width=1\textwidth]{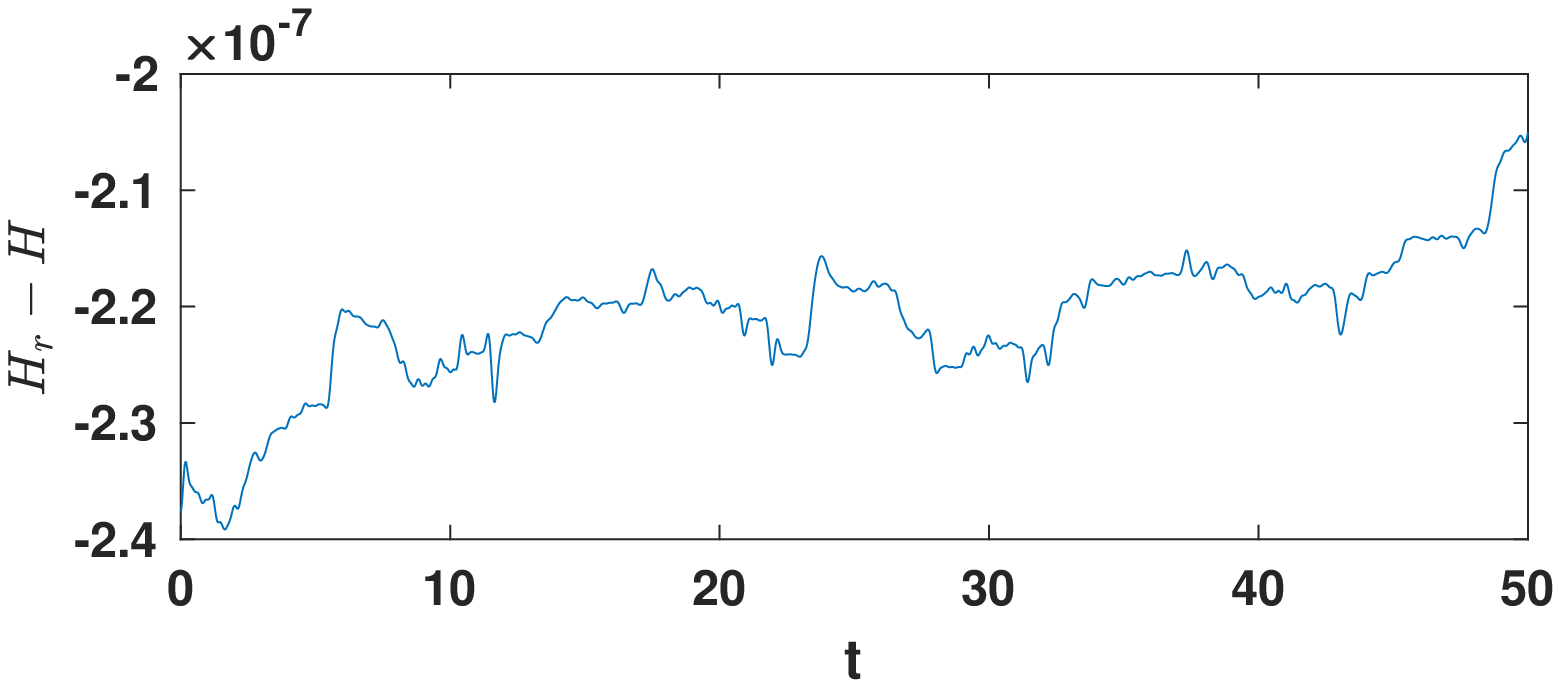}
\end{minipage}
\caption{
$H_r\Delta x-H\Delta x$ of the Standard POD-ROM: $r=10$ (left) and $r=20$ (right). 
}\label{Fig: podg}
\end{figure}

For $r=10$ and $20$, the respective CPU times for online simulations are 0.765~seconds and 0.979~seconds. The associated maximum errors of the reduced approximations are $\mathcal{E}_\infty= 3.291\times 10^{-2}$ and $8.288\times 10^{-3}$, respectively.
The time evolution of Hamiltonian approximation errors for both cases is shown in Figure \ref{Fig: podg}. 
It is observed that the order of Hamiltonian approximation errors varies from $\mathcal{O}(10^{-5})$ when $r=10$ to $\mathcal{O}(10^{-7})$ when $r=20$. 
As $r$ increases, the discrete Hamiltonian becomes more accurate, but it is not time invariant as the standard G-POD ROM is not structure preserving.

\subsection{The SP-POD Models.} 
Next, we consider the structure-preserving ROMs introduced in Section \ref{sec: alg}. Two SP-POD ROMs are applied: the first one uses the standard POD basis, named SP-POD-1; the other uses the POD basis generated from shifted snapshots, named SP-POD-2. 
\paragraph{SP-POD-1.}  This SP-ROM has the following form: 
\begin{equation*}
\left[
\begin{array}{c}
\dot{\ba}\\
\dot{\bb}
\end{array}
\right] =
\left[
\begin{array}{cc}
0 & \bPhi_u^\top \bPhi_v \\
-\bPhi_v^\top \bPhi_u & 0
\end{array}
\right]
\left[
\begin{array}{c}
-\bA_r \ba + \bPhi_u^\intercal \bg(\bPhi_u \ba) \\
\bb
\end{array}
\right],
\end{equation*}
where $\bA_r = \bPhi_u^\top \bA \bPhi_u$. 
The coefficient matrix on the right-hand-side of the system is skew-symmetric, which has the same structure as that of the full-order model.
Thus, we expect a constant discrete Hamiltonian in the reduced-order simulation. 
Using the symplectic midpoint scheme, the discrete system reads: 
\begin{equation}
\left[
\begin{array}{c}
\frac{\ba^{k+1}-\ba^{k}}{\Delta t}\\
\frac{\bb^{k+1}-\bb^{k}}{\Delta t}
\end{array}
\right] =
\left[
\begin{array}{cc}
0 & \bPhi_u^\top \bPhi_v \\
-\bPhi_v^\top \bPhi_u & 0
\end{array}
\right]
\left[
\begin{array}{c}
-\bA_r \frac{\ba^{k+1}+\ba^k}{2} + \bPhi_u^\intercal \bg(\bPhi_u \frac{\ba^{k+1}+\ba^k}{2}) \\
\frac{\bb^{k+1}+\bb^k}{2}
\end{array}
\right]
\end{equation}
with the initial condition $\ba^0 = \bPhi_u^\top \bu_0$ and $\bb^0 = \bPhi_v^\top \bv_0$.

\begin{figure}[htb]
\centering
\begin{minipage}[ht]{0.45\linewidth}
\includegraphics[width=1\textwidth]{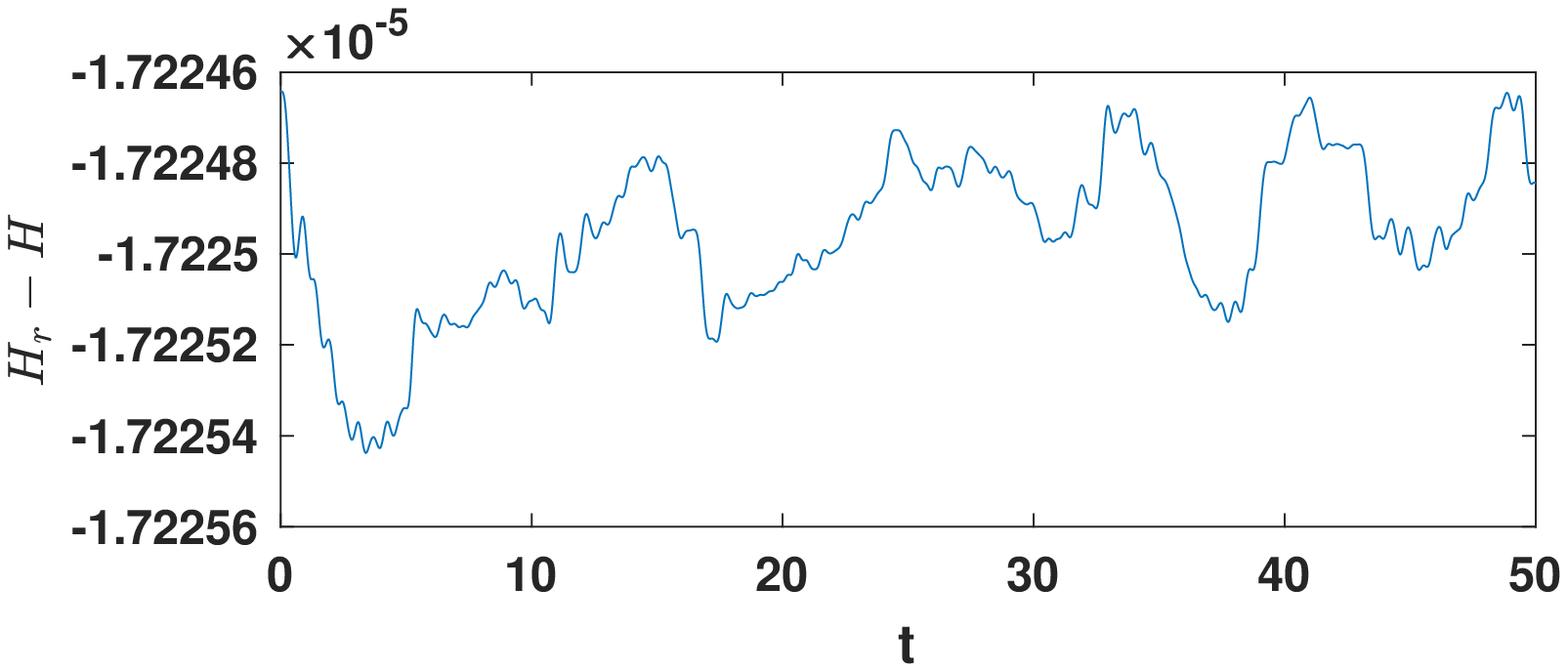}
\end{minipage}
\hspace{1cm}
\begin{minipage}[ht]{0.45\linewidth}
\includegraphics[width=1\textwidth]{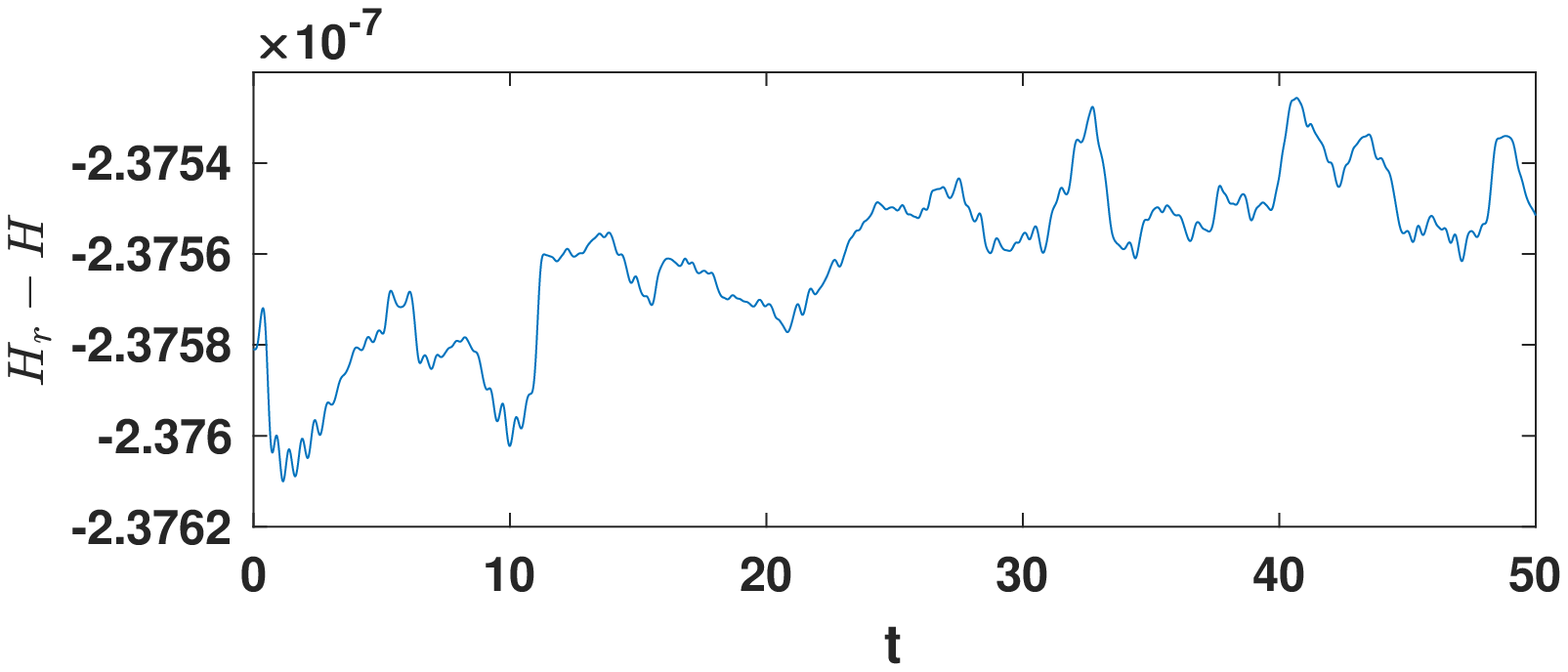}
\end{minipage}
\caption{
$H_r\Delta x-H\Delta x$ of the SP-POD-1: $r=10$ (left) and $r=20$ (right). 
}\label{Fig: sp-pod1}
\end{figure}

It takes $0.601$ seconds and $0.786$ seconds to complete the simulations of SP-POD-1 model with dimensions $r=10$ and $r=20$, respectively. The corresponding maximum errors of the reduced approximations is $\mathcal{E}_\infty= 3.291\times 10^{-2}$ and $8.298\times 10^{-3}$.
Errors of the Hamiltonian approximations in both cases are presented in Figure \ref{Fig: sp-pod1}. 
It is seen that $H_r\Delta x-H\Delta x$ possesses the same magnitude as that in the G-POD model of the same dimensions. 
This is because the POD truncation causes information loss, which results in the Hamiltonian approximation error. Therefore, we next correct the reduced Hamiltonian by using the POD basis from shifted snapshots. 

\paragraph{SP-POD-2.} This SP-ROM has the following form: 
\begin{equation}
\left[
\begin{array}{c}
\dot{\ba}\\
\dot{\bb}
\end{array}
\right] =
\left[
\begin{array}{cc}
0 & \bPhi_u^\top \bPhi_v \\
-\bPhi_v^\top \bPhi_u & 0
\end{array}
\right]
\left[
\begin{array}{c}
-\bA_r \ba -\bPhi_u^\top \bA \bu_0 + \bPhi_u^\intercal \bg(\bPhi_u \ba +\bu_0) \\
\bb + \bPhi_v^\intercal \bv_0
\end{array}
\right].
\end{equation}
Using the same symplectic midpoint method, we have the discrete system as follows. 
\begin{equation}
\left[
\begin{array}{c}
\frac{\ba^{k+1}-\ba^{k}}{\Delta t}\\
\frac{\bb^{k+1}-\bb^{k}}{\Delta t}
\end{array}
\right] =
\left[
\begin{array}{cc}
0 & \bPhi_u^\top \bPhi_v \\
-\bPhi_v^\top \bPhi_u & 0
\end{array}
\right]
\left[
\begin{array}{c}
-\bA_r \frac{\ba^{k+1}+\ba^k}{2} -\bPhi_u^\top \bA \bu_0 + \bPhi_u^\intercal \bg(\bPhi_u \frac{\ba^{k+1}+\ba^k}{2}+\bu_0) \\
\frac{\bb^{k+1}+\bb^k}{2}+ \bPhi_v^\intercal \bv_0
\end{array}
\right]
\end{equation}
with the initial condition $\ba^0 = {\bf 0}$ and $\bb^0 = {\bf 0}$.

\begin{figure}[htb]
\centering
\begin{minipage}[ht]{0.45\linewidth}
\includegraphics[width=1\textwidth]{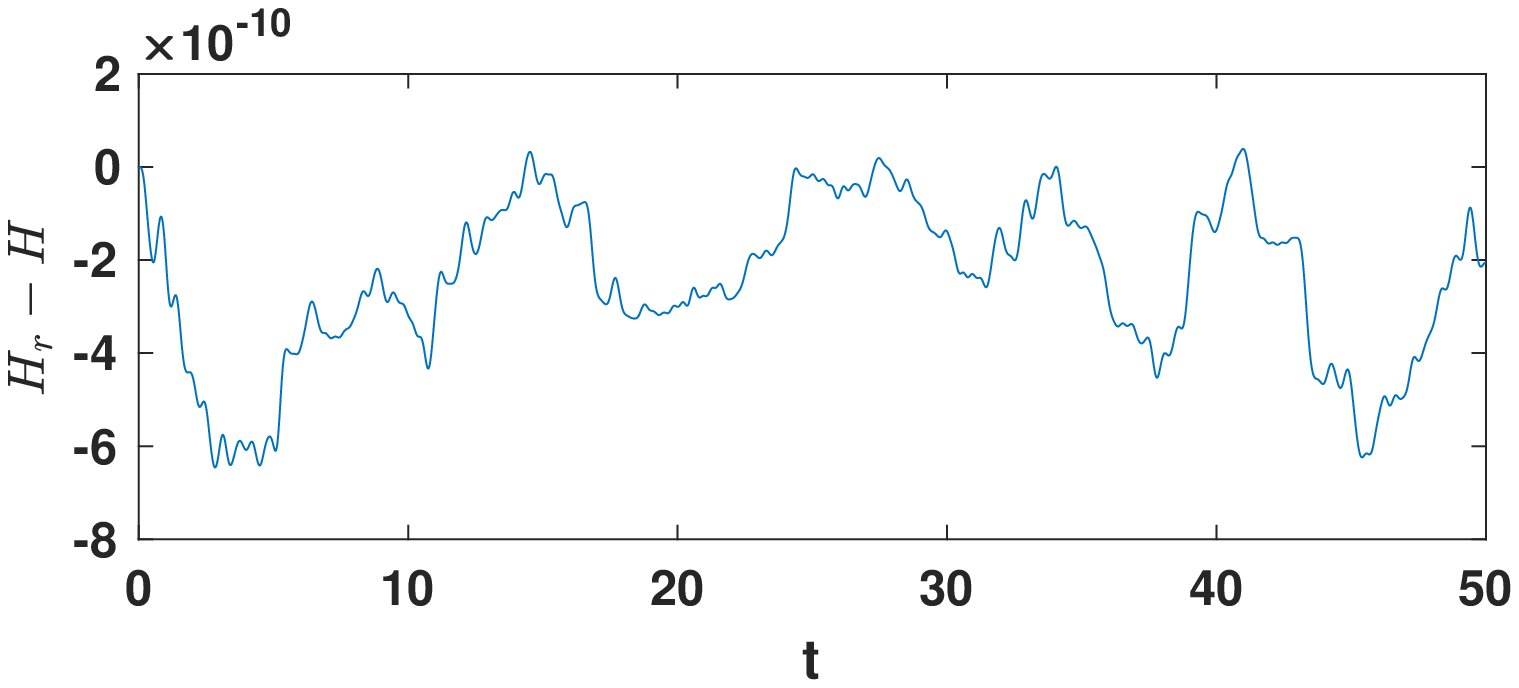}
\end{minipage}
\hspace{1cm}
\begin{minipage}[ht]{0.45\linewidth}
\includegraphics[width=1\textwidth]{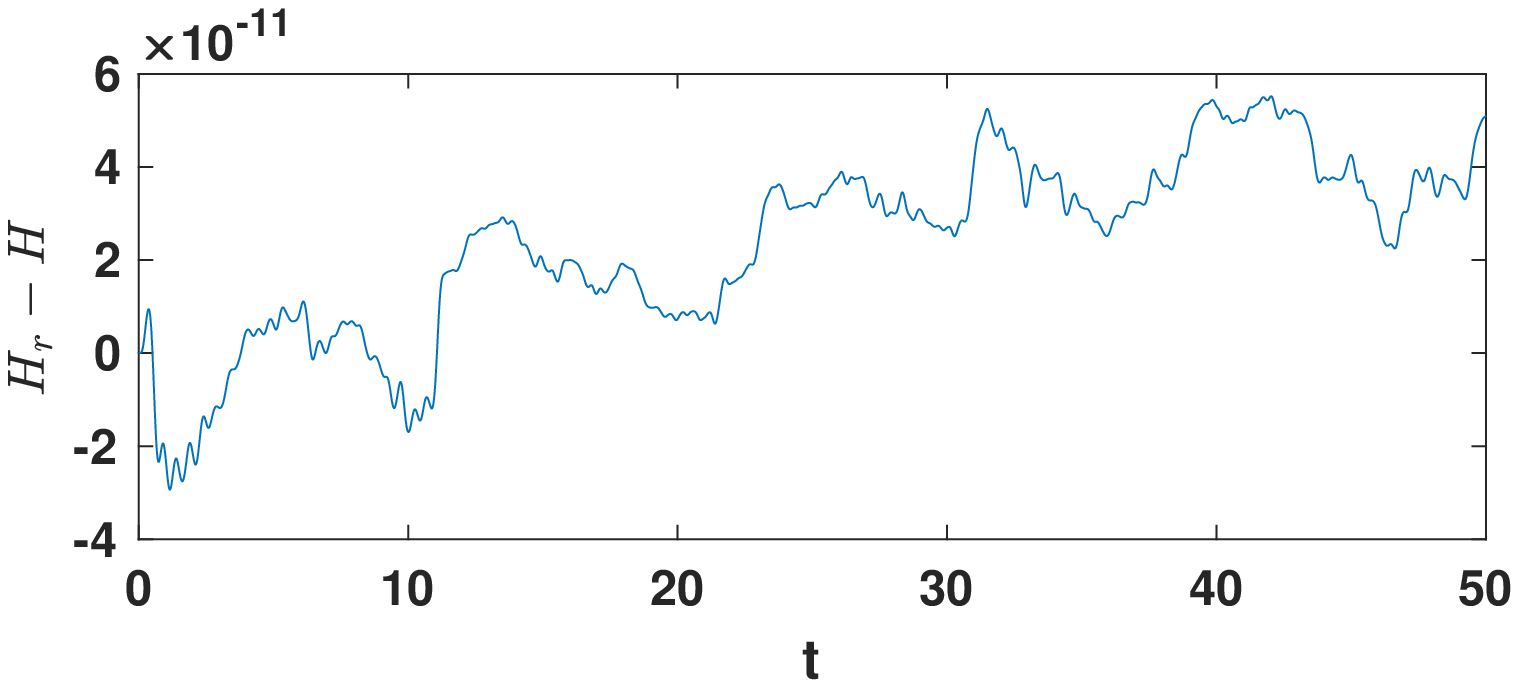}
\end{minipage}
\caption{
$H_r\Delta x-H\Delta x$ of the SP-POD-2: $r=10$ (left) and $r=20$ (right). 
}\label{Fig: sp-pod2}
\end{figure}

The CPU times of the simulations when $r=10$ and $20$ are $0.630$ seconds and $0.805$ seconds, respectively. The corresponding maximum errors are $\mathcal{E}_\infty= 3.711\times 10^{-2}$ and $1.015\times 10^{-2}$.
Figure \ref{Fig: sp-pod2} displays the time history of the discrete Hamiltonian approximation errors during the simulations in both cases. 
Compared with those obtained from SP-POD-1 model, the Hamiltonian approximation errors of SP-POD-2 shrink to $\mathcal{O}(10^{-10})$ and $\mathcal{O}(10^{-11})$ in these two cases. This illustrates the advantage of using POD basis from shifted snapshots for improving the energy approximation. 

\subsection{The SP-DEIM Models\label{sec: SP-DEIM}} The SP-POD models still have the computational complexity dependent on $n$ due to the nonlinearity of $\bg(\cdot)$. 
To reduce the computational cost, we utilize the SP-DEIM ROMs. Two models are applied: the first one is named by SP-DEIM-1 that uses the standard POD and DEIM basis, the other one is named SP-DEIM-2 that uses the POD basis from shifted state snapshots and DEIM basis generated from shifted nonlinear snapshots. 
In both models, we select the number of DEIM basis to be twice as many as that of POD basis. 

\paragraph{SP-DEIM-1.} This model has the following form: 
\begin{equation}
\left[
\begin{array}{c}
\dot{\ba}\\
\dot{\bb}
\end{array}
\right] =
\left[
\begin{array}{cc}
0 & \bPhi_u^\top \bPhi_v \\
-\bPhi_v^\top \bPhi_u & 0
\end{array}
\right]
\left[
\begin{array}{c}
-\bA_r \ba + \bPhi_u^\intercal \bg(\bPhi_u \ba)\mathbb{P}^\intercal C \\
\bb
\end{array}
\right].
\end{equation}
With the symplectic midpoint rule, we have 
\begin{equation}
\left[
\begin{array}{c}
\frac{\ba^{k+1}-\ba^{k}}{\Delta t}\\
\frac{\bb^{k+1}-\bb^{k}}{\Delta t}
\end{array}
\right] =
\left[
\begin{array}{cc}
0 & \bPhi_u^\top \bPhi_v \\
-\bPhi_v^\top \bPhi_u & 0
\end{array}
\right]
\left(
\left[
\begin{array}{c}
-\bA_r \frac{\ba^{k+1}+\ba^k}{2}+ \bPhi_u^\intercal \bg(\bPhi_u \frac{\ba^{k+1}+\ba^k}{2})\mathbb{P}^\intercal C\\
\frac{\bb^{k+1}+\bb^k}{2}
\end{array}
\right]
\right)
\end{equation}
with the initial condition $\ba^0 = \bPhi_u^\intercal \bu_0$ and $\bb^0 = \bPhi_v^\intercal \bv_0$.

\begin{figure}[htb]
\centering
\begin{minipage}[ht]{0.45\linewidth}
\includegraphics[width=1\textwidth]{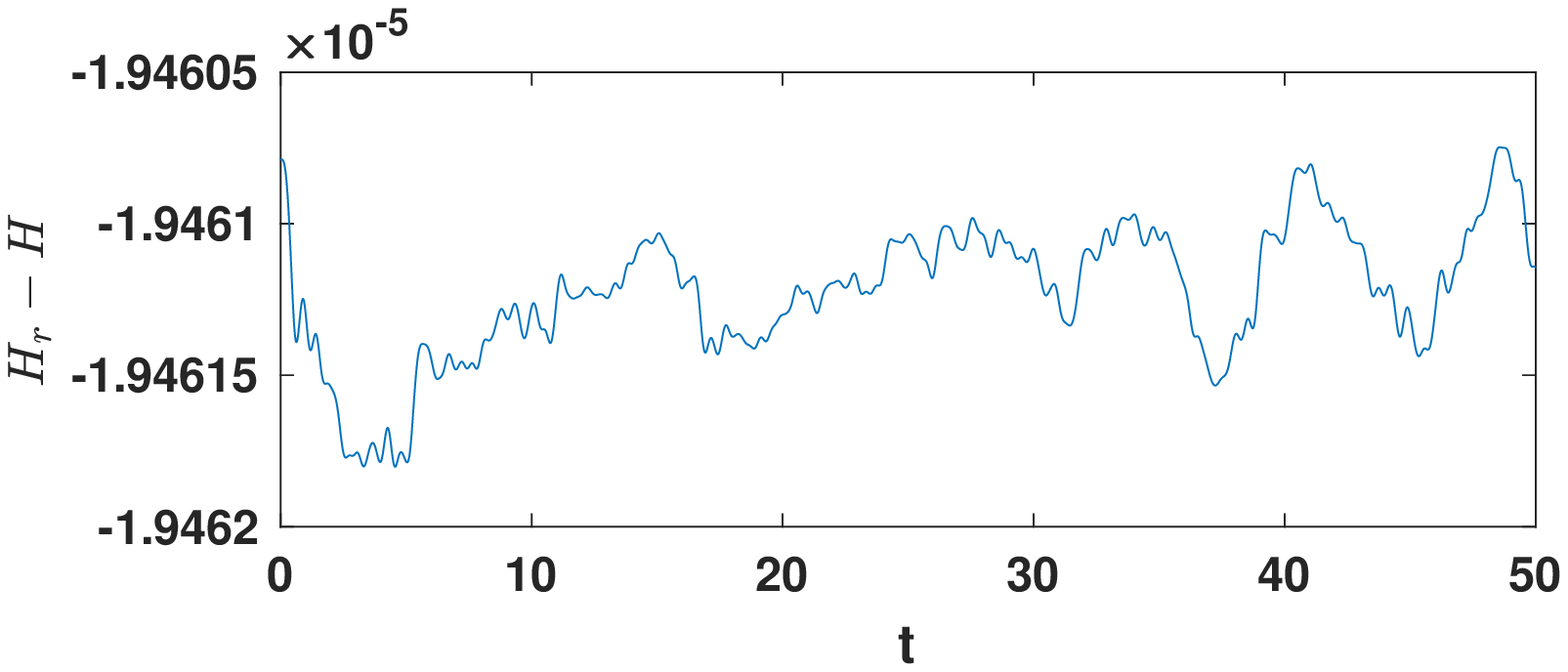}
\end{minipage}
\hspace{1cm}
\begin{minipage}[ht]{0.45\linewidth}
\includegraphics[width=1\textwidth]{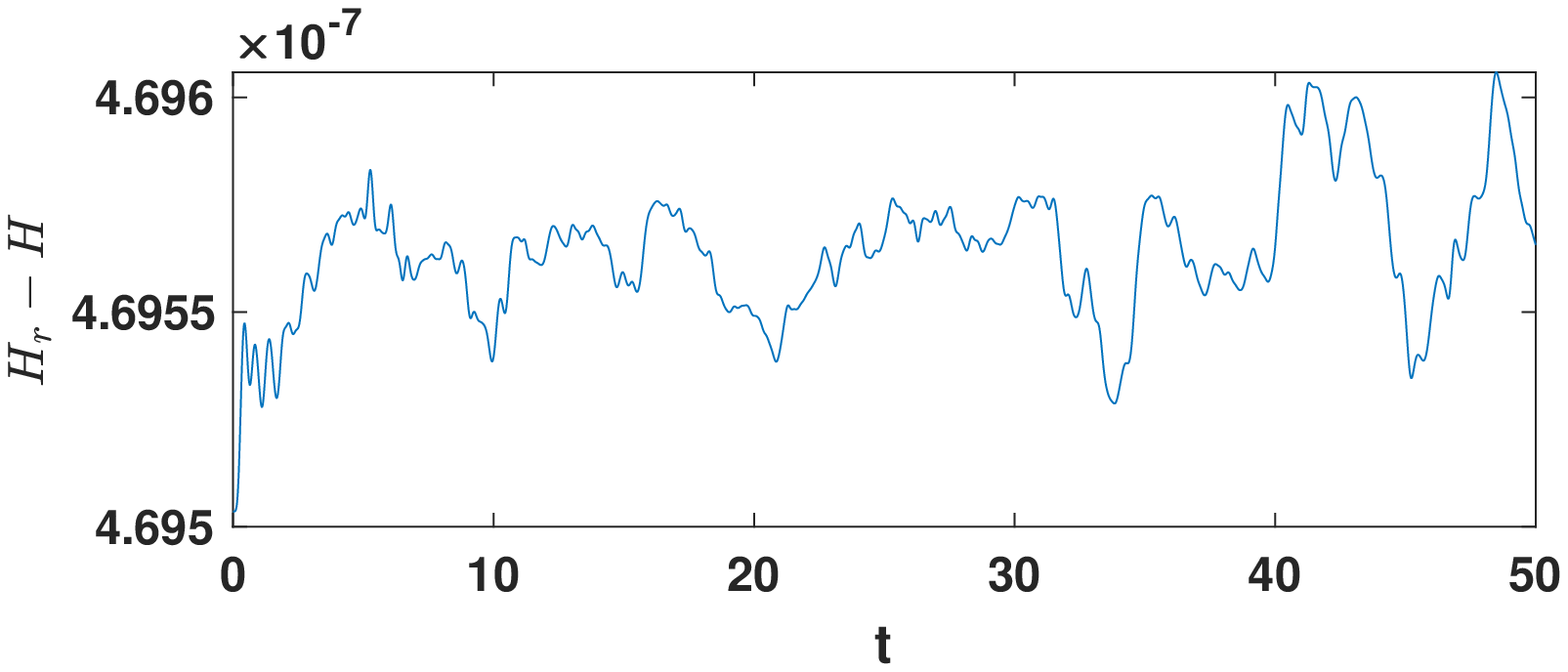}
\end{minipage}
\caption{
$H_r\Delta x - H\Delta x$ in SP-DEIM-1: $r=10$ (left) and $r=20$ (right). 
}\label{Fig: sp-deim1}
\end{figure}

For $r=10$ and $20$, it respectively takes $0.141$ seconds and $0.209$ seconds to complete the simulations. 
Correspondingly, the maximum errors of the reduced-order simulations are $\mathcal{E}_\infty= 3.364\times 10^{-2}$ and $8.473\times 10^{-3}$.
The errors of Hamiltonian function approximations in both cases are presented in Figure \ref{Fig: sp-deim1}. 
Their magnitudes are close to those obtained by the SP-POD-1 model.  

\paragraph{SP-DEIM-2.} To improve the discrete Hamiltonian, we use the POD basis generated from shifted state snapshots and the DEIM basis from shifted nonlinear snapshots. 
The model has the following form:
\begin{equation}
\left[
\begin{array}{c}
\dot{\ba}\\
\dot{\bb}
\end{array}
\right] =
\left[
\begin{array}{cc}
0 & \bPhi_u^\top \bPhi_v \\
-\bPhi_v^\top \bPhi_u & 0
\end{array}
\right]
\left[
\begin{array}{c}
- \bA_r \ba -\bPhi_u^\top \bA\bu_0 + \bPhi_u^\intercal \bg(\bPhi_u \ba + \bu_0)\mathbb{P}^\intercal C \\
\bb + \bPhi_v^\intercal \bv_0
\end{array}
\right].
\end{equation}
After applying the symplectic midpoint rule, we have 
\begin{equation}
\left[
\begin{array}{c}
\frac{\ba^{k+1}-\ba^{k}}{\Delta t}\\
\frac{\bb^{k+1}-\bb^{k}}{\Delta t}
\end{array}
\right] =
\left[
\begin{array}{cc}
0 & \bPhi_u^\top \bPhi_v \\
-\bPhi_v^\top \bPhi_u & 0
\end{array}
\right]
\left[
\begin{array}{c}
-\bA_r \frac{\ba^{k+1}+\ba^k}{2} -\bPhi_u^\top \bA  \bu_0 + \bPhi_u^\intercal \bg(\bPhi_u \frac{\ba^{k+1}+\ba^k}{2} + \bu_0)\mathbb{P}^\intercal C\\
\frac{\bb^{k+1}+\bb^k}{2} + \bPhi_v^\intercal \bv_0
\end{array}
\right]
\end{equation}
with the initial condition $\ba^0 = {\bf 0}$ and $\bb^0 = {\bf 0}$.

\begin{figure}[htb]
\centering
\begin{minipage}[ht]{0.45\linewidth}
\includegraphics[width=1\textwidth]{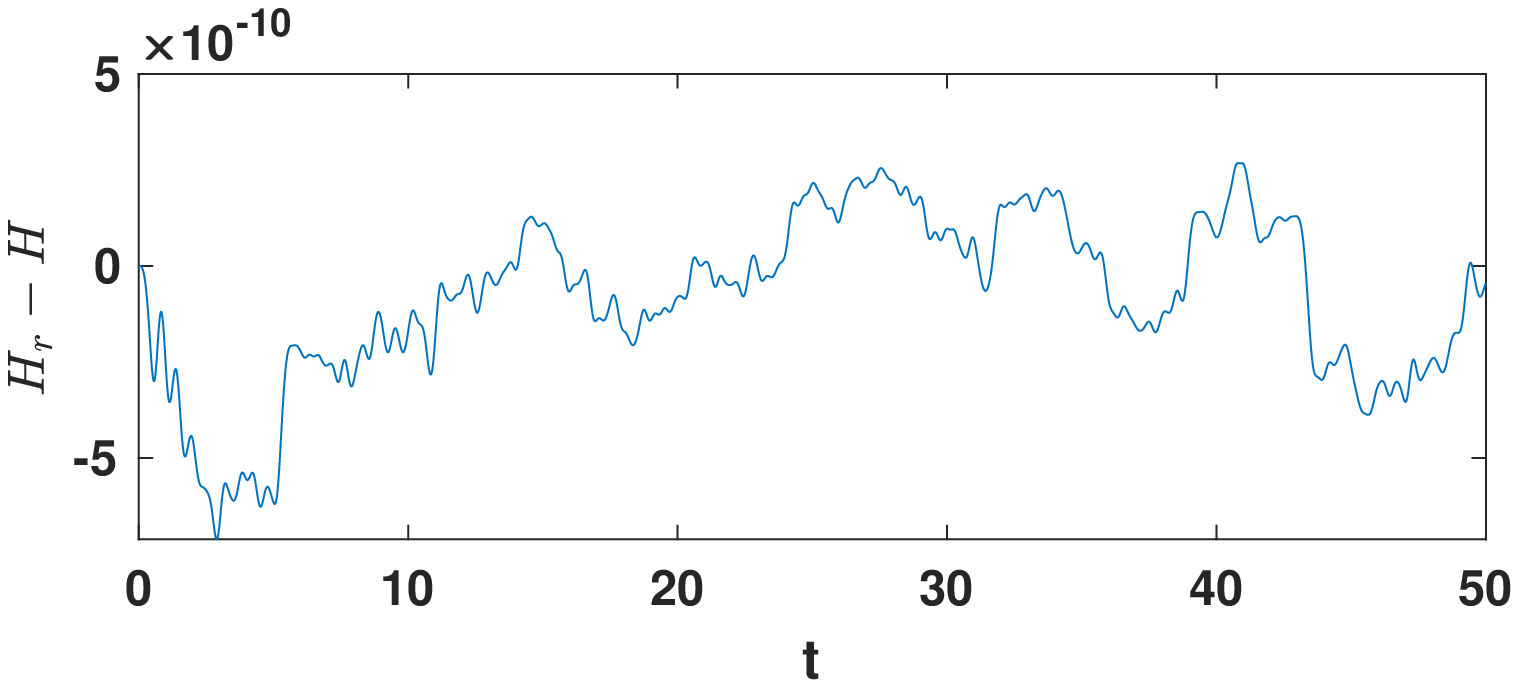}
\end{minipage}
\hspace{1cm}
\begin{minipage}[ht]{0.45\linewidth}
\includegraphics[width=1\textwidth]{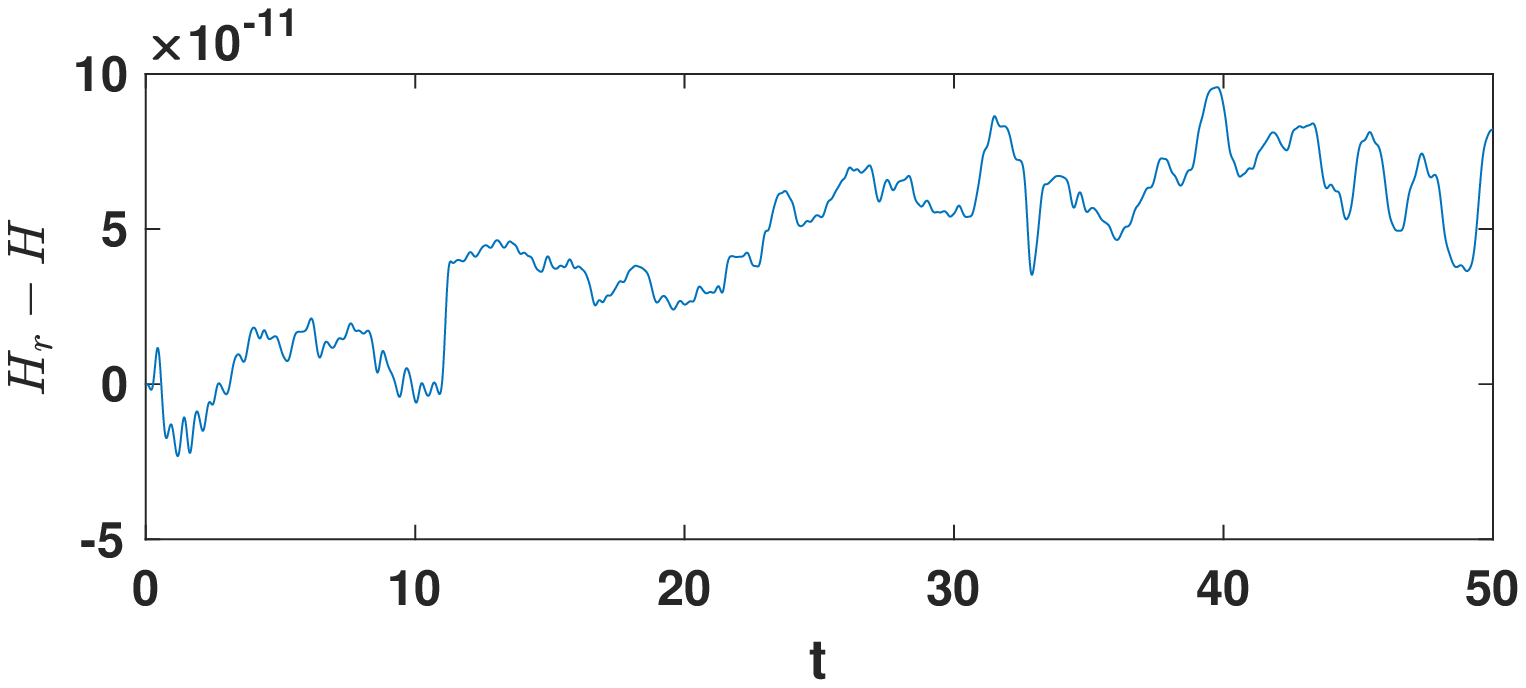}
\end{minipage}
\caption{
$H_r\Delta x-H\Delta x$ in SP-DEIM-2: $r=10$ (left) and $r=20$ (right). 
}\label{Fig: sp-deim2}
\end{figure}

The CPU times for online simulations when $r=10$ and $20$ are $0.143$~seconds and $0.210$~seconds, respectively. The associated maximum errors are $\mathcal{E}_\infty= 3.490\times 10^{-2}$ and $1.311\times 10^{-2}$. 
The errors are bigger than those obtained from the SP-POD-2 model, which can be improved by increasing the number of DEIM basis and interpolation points. 
The Hamiltonian function approximation errors for both cases are presented in Figure \ref{Fig: sp-deim2}, whose magnitude are close to those in the SP-POD-2 simulations. It illustrates the SP-DEIM-2 model is able to preserve the discrete Hamiltonian.

\subsection{Summary of numerical experiments}

We summarize the above test cases in Tables \ref{tab: ex1}-\ref{tab: ex2}. The results are benchmarked by the full order solutions:  
the discrete Hamiltonian value is $H\Delta x = 1.258\times 10^{-1}$ and the CPU time for the full order simulation is $t_{cpu}= 52.8$ seconds.
\begin{table}[htp]
\begin{center}
\caption{Nonlinear wave equations: comparison of ROMs (r=10)}
\label{tab: ex1}
\begin{tabular}{| c | c | c | c | c | c | c |} \hline
  {}                    &   G-ROM                 & SP-POD-1               & SP-POD-2                & SP-DEIM-1             & SP-DEIM-2 \\ \hline
  $\mathcal{E}_\infty$  &  $3.291\times 10^{-2}$  & $3.291\times 10^{-2}$  & $3.711\times 10^{-2}$   & $3.365\times 10^{-2}$ & $3.490\times 10^{-2}$   \\
  $H_r\Delta x-H\Delta x$    			& $\mathcal{O}(10^{-5})$  & $\mathcal{O}(10^{-5})$ & $\mathcal{O}(10^{-10})$ & $\mathcal{O}(10^{-5})$& $\mathcal{O}(10^{-10})$ \\     
  $t_{cpu}$ (s)    		& 0.755                   & 0.601                  & 0.630                   & 0.140                 & 0.143 \\ 
  \hline
  
\end{tabular}
\end{center}
\end{table}

\begin{table}[htp]
\begin{center}
\caption{Nonlinear wave equations: comparison of ROMs (r=20)}
\label{tab: ex2}
\begin{tabular}{| c | c | c | c | c | c |} \hline
  {}                    &   G-ROM                 & SP-POD-1               & SP-POD-2                & SP-DEIM-1             & SP-DEIM-2 \\ \hline
  $\mathcal{E}_\infty$  &  $8.288\times 10^{-3}$  & $8.298\times 10^{-3}$  & $1.152\times 10^{-2}$   & $8.473\times 10^{-3}$ & $1.311\times 10^{-2}$   \\
  $H_r\Delta x-H\Delta x$    			& $\mathcal{O}(10^{-7})$  & $\mathcal{O}(10^{-7})$ & $\mathcal{O}(10^{-11})$ & $\mathcal{O}(10^{-7})$& $\mathcal{O}(10^{-11})$ \\     
  $t_{cpu}$ (s)    		& 0.979                   & 0.786                  & 0.805                   & 0.209                 & 0.210 \\ 
  \hline
  
\end{tabular}
\end{center}
\end{table}

Based on the results, we can draw the following conclusions:
(i) the structure-preserving ROMs with basis from shifted snapshots (SP-POD-2 and SP-DEIM-2)  are able to preserve the energy, which achieve better Hamiltonian approximations than other ROMs;
(ii) using the SP-DEIM models reduces the CPU time of their SP-POD counterparts; The speedup factor could increase if the nonlinear function is more complicated or the dimension of the full order model is bigger; 
(iii) using the SP-DEIM models achieves approximation errors close to their SP-POD counterparts; Such differences would be negligible when more DEIM basis and interpolation points are included; 
and
(iv) overall, SP-DEIM-2 outperforms the other ROMs discussed in this paper in terms of discrete Hamiltonian approximation, the accuracy and efficiency of reduced order simulations.

\section{Conclusions}
\noindent \indent Energy preserving schemes have been developed for simulating Hamiltonian PDEs, for which one significant property is to preserve the Hamiltonian function. When model reduction techniques such as the POD method is applied, the standard Galerkin projection would destroy this property, thus the discrete Hamiltonian is not well preserved.  
In \cite{gong2017structure}, the structure-preserving POD have been developed that leads to a constant Hamiltonian approximation.  
However, the computational complexity of the SP-POD ROMs is still high when the gradient of Hamiltonian has non-polynomial  nonlinearities. 
In this paper, we introduce a new structure-preserving DEIM that improves the online efficiency while obtaining an accurate Hamiltonian approximation from the reduced order simulations. 
Numerical experiments demonstrate the efficacy of the proposed approach.

\section{Acknowledgements}
Z. Wang's research was partially supported by the National Science Foundation through grants DMS-1913073 and DMS-2012469, and by the US Department of Energy through grant SC0020270. 

\bibliographystyle{abbrv} 
\bibliography{bib_structure,bib_thesis}
\end{document}